\documentclass{IEEEtran4PSCC}
\usepackage{amssymb}
\usepackage{amsthm}
\theoremstyle{remark}

\ifCLASSINFOpdf
 \usepackage[pdftex]{graphicx}
\else
\fi
\usepackage[cmex10]{amsmath}
\usepackage{booktabs}
\ifCLASSOPTIONcompsoc
  \usepackage[caption=false,font=normalsize,labelfont=sf,textfont=sf]{subfig}
\else
  \usepackage[caption=false,font=footnotesize]{subfig}
\fi
\begin{document}
\title{Loss Induced Maximum Power Transfer\\ in Distribution Networks}

%% To specify the authors when (number of affiliations <= 2)
\author{
\IEEEauthorblockN{Matthew Deakin, \textit{Student Member, IEEE}, Thomas Morstyn, \textit{Member, IEEE},\\ Dimitra Apostolopoulou, \textit{Member, IEEE}, Malcolm McCulloch, \textit{Senior Member, IEEE}}
\IEEEauthorblockA{Department of Engineering Science, University of Oxford, Oxford, UK}
%{\{matthew.deakin, malcolm.mcculloch\}@eng.ox.ac.uk}
}

\maketitle

\begin{abstract}
In this paper, the power flow solution of the two bus network is used to analytically characterise maximum power transfer limits of distribution networks, when subject to both thermal and voltage constraints. Traditional analytic methods are shown to reach contradictory conclusions on the suitability of reactive power for increasing power transfer. Therefore, a more rigorous analysis is undertaken, yielding two solutions, both fully characterised by losses. The first is the well-known thermal limit. The second we define as the `marginal loss-induced maximum power transfer limit'. This is a point at which the marginal increases in losses are greater than increases in generated power. The solution is parametrised in terms of the ratio of resistive to reactive impedance, and yields the reactive power required. The accuracy and existence of these solutions are investigated using the IEEE 34 bus distribution test feeder, and show good agreement with the two bus approximation. The work has implications for the analysis of reactive power interventions in distribution networks, and for the optimal sizing of distributed generation.
\end{abstract}


\section{Introduction}

The number of large scale, low-carbon generators in power systems has steadily increased in recent years. For example, in the UK, 46 \% of solar photovoltaic capacity is generated at sites of a size $> 5$ MW \cite{decc2017national}. Given geographical constraints, generation may be located in existing distribution networks that are distant from strong transmission networks. If large amounts of generation are connected to a network, then feeder voltage or thermal limits will eventually be reached, yielding a power transfer limit for the given network (without further network interventions). This can result in curtailment of power generated by low-carbon sources such as solar.

Therefore, one particular problem that has been studied is how to increase the maximum real power that can be transferred through a distribution network subject to voltage and thermal constraints. Sinking reactive power is a well-established method of reducing network voltages to increase the real power that can be transferred \cite{stetz2013improved}, but is known to increase losses \cite{turitsyn2011options}. Unfortunately, previous works do not consider this problem analytically \cite{keane2013state}, and so results are generally simulation-based  \cite{stetz2013improved}. This makes analysis of the wide range of distribution network types and generator impacts difficult to infer from (relatively) small numbers of detailed case studies. Furthermore, the analysis of losses (which reduce net real power transfer) tend to be left as an after-thought - their study is instead usually motivated by government regulation \cite{kalambe2014loss}. For example, a recent EPRI communiqu\'{e} on the topic of increasing real power transfer capabilities discusses only that losses may be \textit{reduced} in feeders with distributed generation \cite{epri2015distribution}, without pointing out that the marginal benefit of distributed generation drops off at high penetration levels (i.e., at power transfer limits). The European IGREENGrid project  \cite{Varela2017igreengrid} also advocates the use of reactive power, but losses are barely mentioned (stating that losses increase by 1\% - 10\% with the proposed network solutions).  These gaps motivate our search for an analytic method to calculate the maximum real power that can be transferred through a distribution network, with the impact of reactive power on losses made explicit.

As such, in this work we investigate (analytically) how losses vary as we increase the amount of generation connected to a feeder, for distributed generators with access to \textit{arbitrary} amounts of reactive power. To do so, we study the closed-form power flow solution of a two bus network. This yields our main result: the existence of the `marginal loss-induced maximum power transfer limit', complementing traditional `thermal' maximum power transfer. The existence of these limits are demonstrated on the IEEE 34 distribution test system, and the work holds corollaries on the utility of reactive power provision and curtailment.

\section{Modelling of Radial Distribution Networks and Inadequecy of Heuristic Analytic Analyses}\label{s:modelling}

In this section we first define and solve the two-bus power flow problem. We then use traditional (heuristic) analyses to consider the maximum power transfer problem, demonstrating an inconsistency in the conclusions that these analyses draw.

\begin{figure}
\centering
\includegraphics[width=0.32\textwidth]{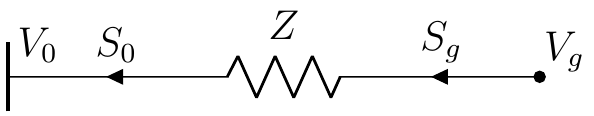}
\caption{Two bus power flow model (sign convention as indicated).}
\label{f:transmission_plot}
\end{figure}

With reference to Figure \ref{f:transmission_plot}, we consider the impedance of the line $Z = R + \text{j}X$ where $R, X$ are both non-negative real numbers. We use $(\cdot )^{*}$ to denote complex conjugate and $|\cdot |$ to denote the magnitude. The $R/X$ ratio is defined as $\lambda = R/X$. We use the notation $S_{(\cdot)} = P_{(\cdot)} + \text{j}Q_{(\cdot)}$, with $P_{(\cdot)},Q_{(\cdot)} \in \mathbb{R} $ real and reactive powers. $V_{0}, S_{0}$ and $V_{g},S_{g}$ representing the (complex) voltage and apparent power at the reference and generator bus respectively.

\subsection{Two bus load flow solution}

Kirchoff's law stipulates $V_{g} = V_{0} - ZI_{0}$, and considering the identity $S_{0} = V_{0}I_{0}^{*}$ we derive the equation

\begin{equation}\label{e:eqn0}
V_{g} = V_{0} - Z \dfrac{S_{0}^{*}}{V_{0}^{*}}\,,
\end{equation}

\noindent where $V_{g}, S_{g} \in \mathbb{C}$, and, without loss of generality, $V_{0} \in \mathbb{R}$. Conservation of energy stipulates that $S_{0} = S_{g} - S_{l}$, and so we can therefore write down the network losses as 
\begin{equation*}
S_{l} = Z\dfrac{|S_{g}|^{2}}{|V_{g}|^{2}}\,. 
\end{equation*}
We now consider a change of co-ordinates,
\begin{equation}\label{e:rot_scale}
\tilde{S}_{(\cdot)} = S_{(\cdot)}Z^{*},
\end{equation}
\noindent which represent a rotation and scaling operation \cite{vournas2015maximum}.  These result in a modified set of equations
\begin{align}
V_{g} =& \, V_{0} - \dfrac{\tilde{S}_{0}^{*}}{V_{0}^{*}}\,, \label{e:deriv1} \\
\tilde{S}_{l} =& \, \dfrac{|\tilde{S}_{g}|^{2}}{|V_{g}|^{2}}\,, \label{e:deriv2}
\end{align}
\noindent with solution in voltage and losses (proportional to current) as
\begin{align}
|V_{g}|^{2} &= \tilde{P}_{g} + \dfrac{V_{0}^{2}}{2} \pm \sqrt{\dfrac{V_{0}^{4}}{4} + V_{0}^{2}\tilde{P}_{g} - \tilde{Q}_{g}^{2}},\label{e:v2}\\
\tilde{S}_{l} &= \tilde{P}_{g} + \dfrac{V_{0}^{2}}{2} \mp \sqrt{\dfrac{V_{0}^{4}}{4} + V_{0}^{2}\tilde{P}_{g} - \tilde{Q}_{g}^{2}}\label{e:sl}.
\end{align}
%\begin{eqnarray}
%|V_{g}|^{2} &= \tilde{P}_{g} + \dfrac{V_{0}^{2}}{2} \pm \sqrt{\dfrac{V_{0}^{4}}{4} + V_{0}^{2}\tilde{P}_{g} - \tilde{Q}_{g}^{2}},\label{e:v2}\\
%\tilde{S}_{l} &= \tilde{P}_{g} + \dfrac{V_{0}^{2}}{2} \mp \sqrt{\dfrac{V_{0}^{4}}{4} + V_{0}^{2}\tilde{P}_{g} - \tilde{Q}_{g}^{2}}\label{e:sl}.
%\end{eqnarray}
\noindent We note that the solution to \eqref{e:v2} and \eqref{e:sl} are non-negative real numbers (there does not exist a solution to the power flow equations if the discriminants are not positive). It is well known that the power flow equations has multiple solutions (see, e.g., \cite{hiskens1995analysis}): here we see in the two bus case that this results in a `high voltage, low loss' solution, and a `low voltage, high loss' solution. Under nominal conditions it is self-evident that we should operate in the former. In either case (i.e. $\pm \to +$ or $\pm \to -$ in \eqref{e:v2} and \eqref{e:sl}) we note that the following identity holds:

\begin{equation}
\tilde{S}_{l} + |V_{g}|^{2} - V_{0}^{2} - 2\tilde{P}_{g} = 0 \, . \label{e:id_1}
\end{equation}

%\begin{figure}
%\centering
%\includegraphics[width=0.37\textwidth]{vg_slt.pdf}
%\caption{High voltage, low loss solution to the two bus voltage case in normalised co-ordinates $\tilde{S}_{g} = S_{g}Z^{*}$.}\label{f:vg_slt}
%\end{figure}

\subsection{Inadequecy of Maximum Power Transfer Theorem}\label{s:limitations}

The Maximum Power Transfer Theorem (MPTT) is a well-known solution to the problem of maximum power that can be transferred to a load from a voltage source through an impedance \cite{mclaughlin2007deglorifying}. The problem definition only considers the existence of a maximum loading point and not on any operating constraints. Although the limitations of the MPTT in power systems are well known \cite{mclaughlin2007deglorifying}, it might be assumed that this might provide a bound on the power $P_{0}$ that can be transferred. In which case, the problem can be stated as
\begin{align*}
\max _{\tilde{S}_{g}} & \quad P_{0}\\
\text{s.t. } & \quad \dfrac{V_{0}^{4}}{4} + V_{0}^{2}\tilde{P}_{g} - \tilde{Q}_{g}^{2} \, \geq \, 0.
\end{align*}
%\begin{eqnarray*}
%\max _{\tilde{S}_{g}} &P_{0}\\
%\text{s.t. } &\dfrac{V_{0}^{4}}{4} + V_{0}^{2}\tilde{P}_{g} - \tilde{Q}_{g}^{2} \geq 0.
%\end{eqnarray*}
\noindent Using \eqref{e:rot_scale}, we can write $P_{0}$ as
\begin{align*}
P_{0} 	=& \, \dfrac{1}{|Z|^{2}}(R\tilde{P}_{0} - X\tilde{Q}_{0}),\\
		=& \, \dfrac{1}{|Z|^{2}}\big(R(\tilde{P}_{g} - \tilde{S}_{l}) - X\tilde{Q}_{g}\big),
\end{align*}
\noindent and finally by operating in a `low losses/high voltage' region (such that $\mp \to -$ in \eqref{e:sl}) that
\begin{equation*}
P_{0} =  \dfrac{R}{|Z|^{2}}\Big( -\dfrac{V_{0}^{2}}{2}  + \sqrt{\dfrac{V_{0}^{4}}{4} + V_{0}^{2}\tilde{P}_{g} - \tilde{Q}_{g}^{2}} \Big) - \dfrac{X\tilde{Q}_{g}}{|Z|^{2}}.
\end{equation*}
\noindent By setting $\tilde{Q}_{g} = 0$, we can increase $\tilde{P}_{g}$ and make $P_{0}$ arbitrarily large. Therefore, without operating constraints, the maximum power that can be sent from a generator to a strong grid is unbounded.

Therefore, we can conclude that operational constraints will always be required to consider the maximum power that can be transferred from a generator to a strong network. This motivates out search for an alternative maximum power transfer criteria.

\subsection{Two Heuristic Voltage Regulation Strategies}\label{s:two_heuristics}

We first consider the example of reactive power control of a network constrained by voltage, with the goal of maximising the power transferred to the grid. 

We note that from \eqref{e:rot_scale} and \eqref{e:id_1} we can calculate the net real power transferred as
\begin{align}
P_{0} =& \, P_{g} - P_{l},\\
		=& \, P_{g} - \dfrac{R}{|Z|^{2}}\tilde{S}_{l}, \label{e:p0_in_sl}\\
 =& \, P_{g} - \dfrac{R}{|Z|^{2}}(V_{0}^{2} + 2\tilde{P}_{g} - |V_{g}|^{2}),
\end{align}
\noindent finally yielding
\begin{equation}\label{e:p0_id}
\begin{split}
P_{0} = P_{g}\Big( \dfrac{\lambda ^{-1} - \lambda }{\lambda  + \lambda ^{-1}} \Big) - Q_{g}&\dfrac{2}{\lambda  + \lambda ^{-1}}  \\
& + \dfrac{(V_{0}^{2} + |V_{g}|^{2})}{|Z|}\dfrac{\lambda}{\sqrt{\lambda^{2}+1}}.
\end{split}
\end{equation}
\noindent As $P_{0}$ is here parametrised in terms of $\lambda$, we can study the solution for mostly inductive or resistive lines. As we shall see, these tend towards operational modes that are traditionally used in networks in different settings.

\subsubsection{Unity Power Factor (`UPF') Control} 

Consider a resistive line such that $\lambda  \to \infty$. We see that
\begin{equation}
\lim _{\lambda \to \infty} P_{0} = - P_{g} + \dfrac{V_{0}^{2} + |V_{g}|^{2}}{|Z|}\,.\label{e:upf_control}
\end{equation}
%\begin{align}
%\lim _{\lambda \to \infty} P_{0} =& \, P_{g}\Big( \dfrac{- \lambda }{\lambda} \Big) - Q_{g}\dfrac{2}{\lambda}  + \dfrac{(V_{0}^{2} + |V_{g}|^{2})}{|Z|}\dfrac{\lambda}{\sqrt{\lambda^{2}}}\\
%								 =& \, - P_{g} + \dfrac{V_{0}^{2} + |V_{g}|^{2}}{|Z|}\label{e:upf_control}
%\end{align}
\noindent From \eqref{e:upf_control} we see that if we are at the voltage limit (such that $(V_{0}^{2} + |V_{g}|^{2})/|Z|$ is a positive constant), then we need to \textit{minimise} $P_{g}$ along the locus of points described by \eqref{e:v2} (all of which give a constant value of $|V_{g}|$). In other words, even if the marginal cost of power is zero, we should refrain from using reactive power to generate additional power, because the real power losses increase at a greater rate than the rate at which we can generate additional real power. 

Using little or no reactive power agrees with how distribution networks are traditionally operated (until recently IEEE 1547 stipulated that distributed generators could not participate in volt/var control \cite{farivar2012optimal}). We refer to this heuristic of setting $Q_{g} = 0$ as `UPF' control. That is, once a voltage limit is reached, additional generated real power is curtailed. 

\subsubsection{Solution Boundary (`Bdry') Control}

In the case $\lambda  \to 0$ we can repeat this analysis. Equation \eqref{e:p0_id} becomes
\begin{equation*}
\lim _{\lambda \to 0} P_{0} = P_{g}\,.
\end{equation*}
\noindent Therefore, reactive power should be used to simply retain the voltage within limits, as there is no cost (in a real power sense) of increasing the reactive power. 

If we continue to increase the real power, we will however eventually reach the stability boundary (i.e. `critical points', or the `knee' of the $P-V$ curve \cite{cutsem2007voltage}), which one might assume limits the real power transfer. This occurs when the discriminant in \eqref{e:v2} and \eqref{e:sl} is identically equal to zero (as in the case of the maximum power transfer theorem). In this case, we can therefore use \eqref{e:v2} to derive
\begin{equation}\label{e:p0_stb}
P_{0}^{stb} = \dfrac{V_{0}^{2}}{|Z|}\bigg( - \dfrac{R}{2|Z|} + \dfrac{X}{|Z|}\sqrt{\dfrac{|V_{g}|^{2}}{V_{0}^{2}} - \dfrac{1}{4}} \bigg).
\end{equation}
We note that sometimes it is assumed that this yields a solution which is in some sense `unstable'. In the case of rotor angle stability, this is indeed the case \cite{grainger1994power}; however, in the absence of a model of the generation connected to the grid we cannot explicitly make a judgement on the stability (indeed, stability is inherently a dynamic problem \cite{cutsem2007voltage}).

%, i.e. on the curve
%\begin{equation}\label{e:sln_boundary}
%\dfrac{V_{0}^{4}}{4} + V_{0}^{2}\tilde{P}_{g} - \tilde{Q}_{g}^{2} = 0.
%\end{equation}
%If this is the case we can then write down
%\begin{eqnarray*}
%\tilde{P}_{g} =& |V_{g}|^{2} - \dfrac{V_{0}^{2}}{2},\\
%\tilde{Q}_{g} =& \pm V_{0}\sqrt{|V_{g}|^{2} - \dfrac{V_{0}^{2}}{4}},\\
%S_{l} =& \dfrac{|V_{g}|^{2}}{Z^{*}},
%\end{eqnarray*}
%\noindent yielding

\section{Loss Induced Maximum Power Transfer}\label{s:maximum}

In the previous section we demonstrated that different assumptions about the $R/X$ ratio of the network lead to contradictory advice considering the operation of distributed generation. Therefore, a more general method is required to unite these apparently divergent heuristics. 

\subsection{Problem Statement}

Consider the two bus model shown in Figure \ref{f:transmission_plot} and consider the optimization problem
\begin{subequations}
\label{e:main_1}
\begin{align}
	\max _{S_{g}} & \quad P_{0} \label{e:main_1a}\\ 
	\text{s.t. } & \quad \sqrt{\dfrac{S_{l}}{Z}} \leq I_{+} \label{e:main_1c}\\
				 & \quad |V_{g}| \leq V_{+} \, .\label{e:main_1b}
\end{align}
\end{subequations}
We are hence looking to maximise the maximum power transferred across the network \eqref{e:main_1a}. This is subject to a limit on current in the line \eqref{e:main_1c}. The network is subject to a voltage constraint at the generator \eqref{e:main_1b}. In the sequel we shall assume that the thermal limits are large enough that the voltage limit is encountered first - i.e., that we have long network lines, and thus that at the thermal limits that we are on the locus of points satisfying
\begin{equation}
|V_{g}| = V_{+}\,. \label{e:PVlocus}
\end{equation}
Note also we do not consider bounds on reactive power.

\subsection{Problem Solution}
To solve \eqref{e:main_1}, we first look to maximise the power that can be generated, subject to the thermal limits. 

\newtheorem{thermal_result}{Lemma}
\begin{thermal_result}\label{th:thermal}
The maximum power that can be generated, $\hat{P}_{g}$ , subject to \eqref{e:main_1c}, is defined by 
\begin{equation}
\hat{S}_{g} = \dfrac{\tilde{\hat{S}}_{g}}{Z^{*}}\,,\label{e:lemma1}
\end{equation}
\noindent where
\begin{align}
\tilde{\hat{P_{g}}} =& \, \dfrac{1}{2}(V_{+}^{2} - V_{0}^{2} + |Z|^{2}I_{+}^{2})\,, \label{e:lemma1a}\\
\tilde{\hat{Q}}_{g} =& \, - \sqrt{(V_{+}I_{+}|Z|)^{2} - \tilde{\hat{P_{g}}}^{2}}\,. \label{e:lemma1b}
\end{align}
The real power transferred is 
\begin{equation}
\hat{P}_{0} = \hat{P}_{g} - I_{+}^{2}R \,. \label{e:lemma1c}
\end{equation}
\textit{Proof:} Equation \eqref{e:lemma1a} comes directly from \eqref{e:id_1} and \eqref{e:main_1c}; \eqref{e:lemma1b} comes from substituting this result into \eqref{e:deriv2}, \eqref{e:lemma1} comes by re-substituting into \eqref{e:rot_scale}, and finally, \eqref{e:lemma1c} by conservation of energy. $\square$
\end{thermal_result}

This does not yet solve our optimization problem; this represents only the intersection of \eqref{e:v2} and \eqref{e:sl}. Furthermore, Lemma \ref{th:thermal} has maximised over $P_{g}$ rather than $P_{0}$, and so does not yet resolve the contradictions of the previous section (we shall now show that maximising $P_{g}$ does not necessarily maximise $P_{0}$). We thus turn to our main result.

\newtheorem{main_result}{Theorem}
\begin{main_result}\label{th:main_result}
The solution to \eqref{e:main_1} is at
\begin{equation}
P_{g} = \min \{ \hat{P}_{g}, P_{g}' \}\,,
\end{equation}
where 
\begin{align}
S_{g}' =& \, P_{g}' + \text{j}Q_{g}' = \, \dfrac{\tilde{P}_{g}' + \text{j}\tilde{Q}_{g}'}{Z^{*}}\,,\\
\tilde{P}_{g}' =& \, V_{+}\Big(V_{+} - V_{0}\dfrac{\lambda}{\sqrt{1 + \lambda^{2}}} \Big)\,, \label{e:pgtpr} \\
\tilde{Q}_{g}' =& \, - V_{0}V_{+}\dfrac{1}{\sqrt{1 + \lambda^{2}}}\,,
\end{align}
\noindent and $\hat{P}_{g}$ is as defined in Lemma \ref{th:thermal}. The real power transferred at $S_{g}'$ is 
\begin{equation}
P_{0}' = \dfrac{V_{0}^{2}}{|Z|}\Bigg(\dfrac{V_{+}}{V_{0}} - \dfrac{\lambda}{\sqrt{1 + \lambda^{2}}}\Bigg)\,.\label{e:p0pr}
\end{equation}
\textit{Proof:} See Appendix. $\square$
\end{main_result}

Theorem \ref{th:main_result} describes two ways in which losses can result in a bound on the maximum power that can be transferred. The first, $\hat{P}_{g}$, we refer to as the `thermal' loss-induced limit, as it is induced by thermal limitations imposed by equipment. We now refer to $P_{g}'$ as the `marginal' loss-induced limit. The maximum of \eqref{e:main_1a} is reached when the marginal increase in losses is greater than the marginal increase in generated power $P_{g}$. In the case that $P_{g}' < \hat{P}_{g}$, we can thus conclude that, even if the marginal cost of generating real and reactive power is zero (as can be approximated for the case of solar PV or wind), real power generation should be curtailed beyond $P_{g}'$.

\subsubsection{Comparison with UPF and Solution Boundary Control}

\begin{figure}
\centering
\includegraphics[width=0.42\textwidth]{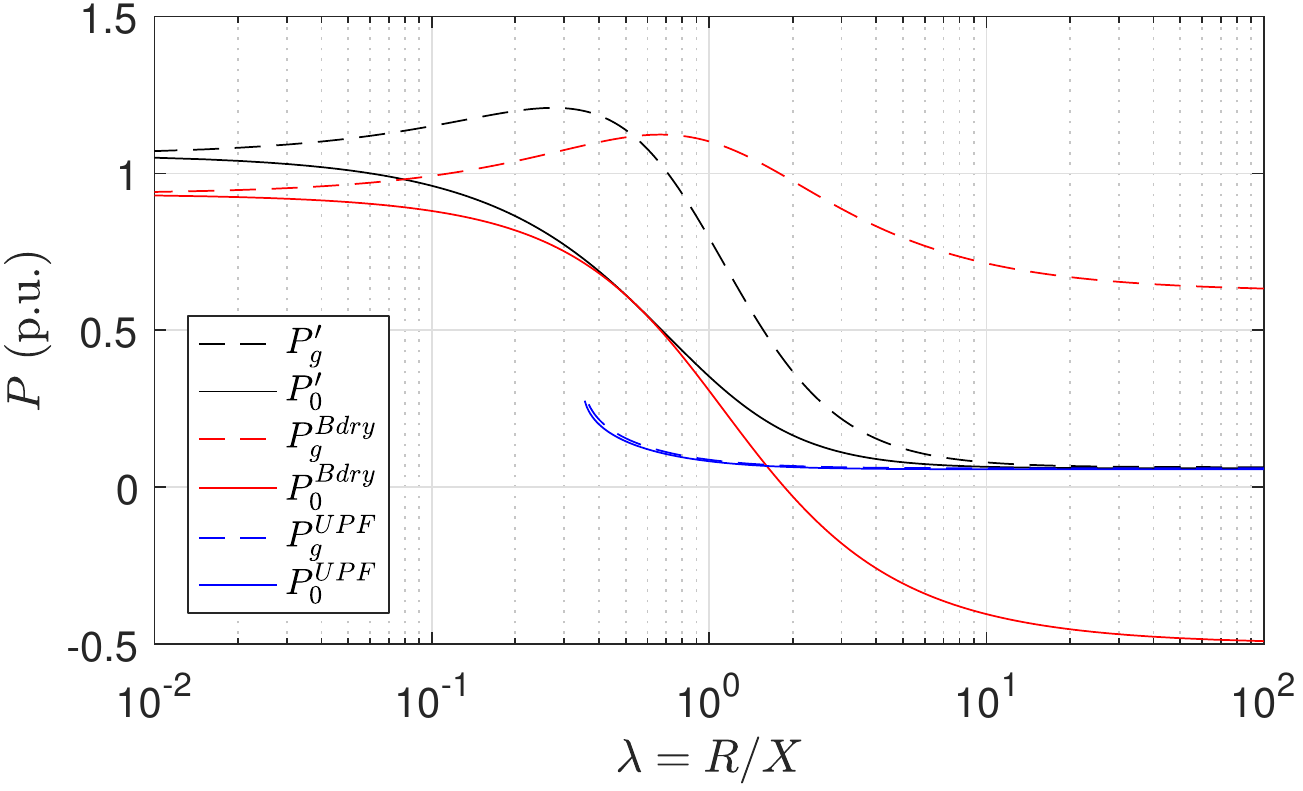}
\caption{Comparison of heuristic (`Bdry', `UPF') methods and the marginal loss-induced maximum power $P_{g}'$ (see Theorem \ref{th:main_result}). Parameters are set as $V_{0} = 1.0, V_{+} = 1.06, |Z| = 1$ (the short circuit current is thus 1 pu). }\label{f:mgtt_1dcomparison}
\end{figure}

For a given $V_{0}, V_{+}, |Z|$ we can compare $P_{g}'$ with the operating strategies considered in Section \ref{s:two_heuristics} (see Figure \ref{f:mgtt_1dcomparison}). We see that for large $\lambda$ the optimal operation approaches that of UPF control, with little use for reactive power. As $\lambda $ is reduced, $P_{g}'$ approaches the solution boundary power $P_{g}^{\text{Bdry}}$. Eventually, for some $\lambda$ we see that $P_{g}' > P_{g}^{\text{Bdry}}$. That is, we have moved from the `high voltage/low losses' region to a `low voltage/high losses' solution (in \eqref{e:v2}, \eqref{e:sl}).

We can calculate the the point at which this happens using \eqref{e:id_1} and \eqref{e:pgtpr}, which yields
\begin{equation}
\lambda' = \dfrac{V_{0}}{\sqrt{4V_{+}^{2} - V_{0}^{2}}}.
\end{equation}
\noindent For the parameters chosen here, $\lambda' = 0.51$, as in Figure \ref{f:mgtt_1dcomparison}. If $V_{0},V_{+}\in [0.9,1.1]$, then $\lambda' \in [0.45,0.77]$. 

As discussed previously, in this region where $P_{g}' > P_{g}^{\text{Bdry}}$, the network is being operated in a region that is traditionally associated with stability issues (either voltage stability in the case of loads, or angle stability in the case of synchronous machines). Although it has not been explicitly been considered in this work, it would be interesting to consider if power electronic interfaced generator could be controlled to overcome this observed `stability' limit (see e.g. \cite{burchett2017Voltage}). We do note that the magnitude of the power sent down the line approaches that of the short circuit power of the line (in this case 1 pu), and so protection issues might need special attention in these cases.

\subsubsection{Thermal efficiency and Power Factor of Marginal Loss-Induced Power Transfer}

Theorem \ref{th:main_result} allows us to study the thermal efficiency and reactive power flows of the network, and allows us to consider the practicalities of operating at this particular point. In particular, thermal efficiency and power factor are two indices that might represent operation that is efficient in some sense. 

\begin{figure}
\centering
\subfloat[]{\includegraphics[width=0.42\textwidth]{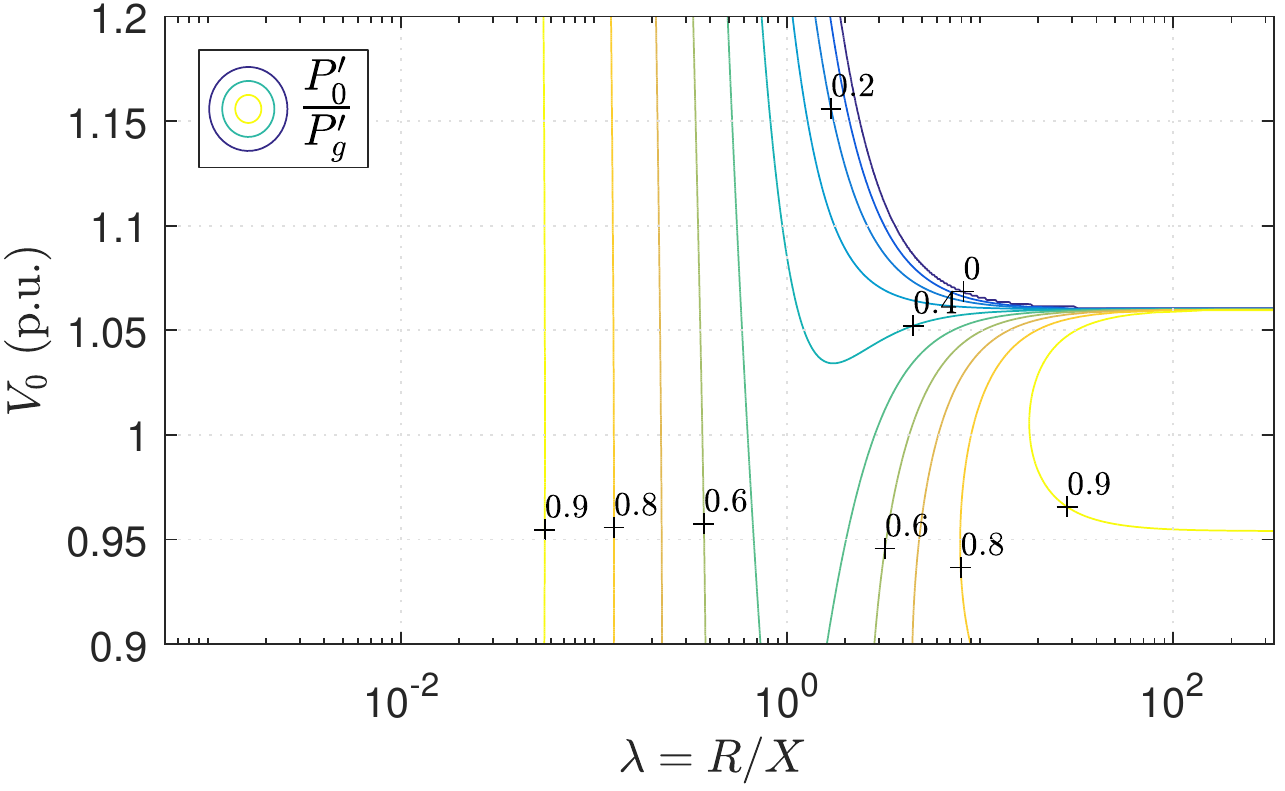}
\label{f:mgtt_eta}}
\hfill
\subfloat[]{\includegraphics[width=0.42\textwidth]{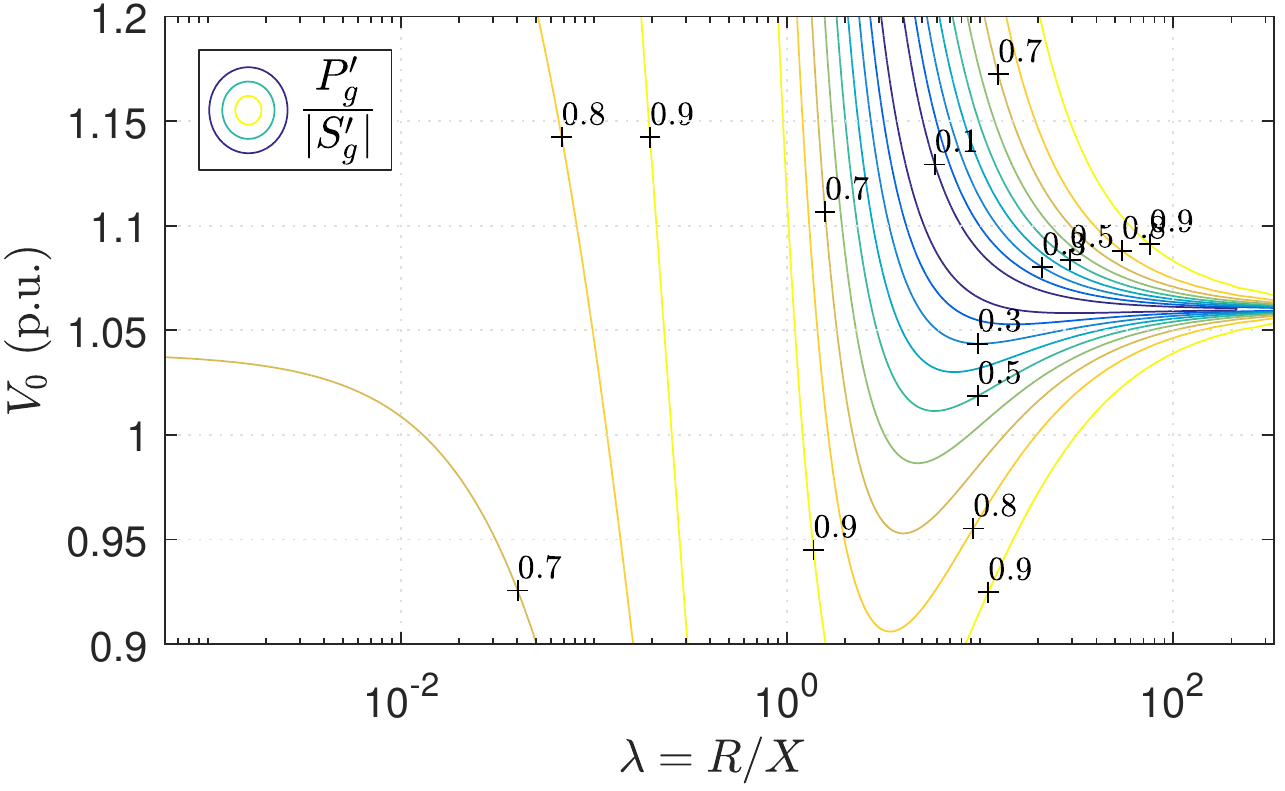}
\label{f:mgtt_pf_g}}
\hfill
\subfloat[]{\includegraphics[width=0.42\textwidth]{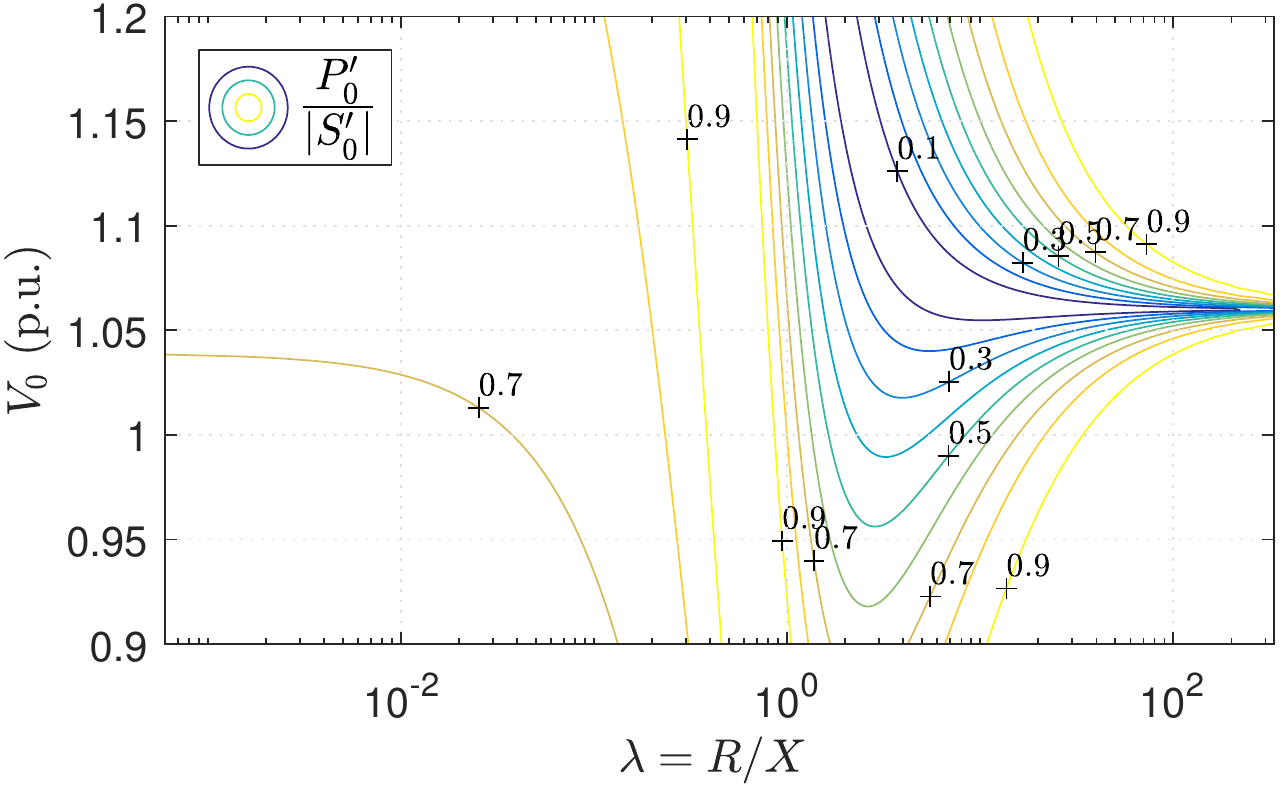}
\label{f:mgtt_pf_0}}
\caption{Performance metrics at the marginal loss induced maximum power transfer limit, with $|Z| = 1, V_{+}=1.06$. (a) Thermal efficiency (zero if positive real power cannot be transferred). (b) Power factor of the generator. (c) Power factor at the substation.}
\end{figure}

The thermal efficiency for given parameters is shown in Figure \ref{f:mgtt_eta}. We see that, for very resistive or very inductive lines, $P_{0}'$ is generally efficient with $P_{0}'/P_{g}' > 90$ \%. However, for moderate values of $\lambda$ the efficiency of power transfer is relatively low, and so unless the marginal cost of power is low, it is unlikely that a line would be operated to this point.

The power factor of the generator (Fig. \ref{f:mgtt_pf_g}) gives a measure of the amount of reactive power that must be supplied to reach this point, either by over-sizing and inverter or by providing shunt elements (switched reactors/capacitors). The power factor at the feeder head is shown in Figure \ref{f:mgtt_pf_0}. Poor power factors at the substation increase the reactive power drawn from the transmission network, which is typically very lossy with respect to $Q$ (this is a strong reason for maintaining substation power factors within bounds). Indeed, it would likely be most efficient to provide reactive power compensation directly at the substation in these cases.

Finally, we also note that there are other metrics and constraints that must be adhered to - the best known are probably protection, harmonic emission standards, and stability issues (as well as lower voltage limits). Furthermore, the thermal and marginal loss-induced limits described here might be limited by the amount of reactive power that is economical to install in a network for intermittent, low-carbon sources such as wind and solar. As such, the reactive power would then become the limiting factor. However, in the event that sufficient reactive power is installed, Theorem \ref{th:main_result} does describe when generated power should be curtailed (even in the case of zero marginal cost power), and bounds the real and reactive power flows that should ever be seen `upstream' in the transmission network.

\section{Distribution Systems Analysis}

As has been previously noted, a majority of PV is large scale in some regions \cite{decc2017national}, and we assume that this trend might continue due to economies of scale. Therefore, we use the preceding analysis to consider the impact on the steady state behaviour of a single, large-scale generator at a single bus of an existing distribution system. In order to simplify analysis, we fix the taps on any in-line voltage regulators.

\begin{figure}
\hspace{2.5em}
\includegraphics[width=0.4\textwidth]{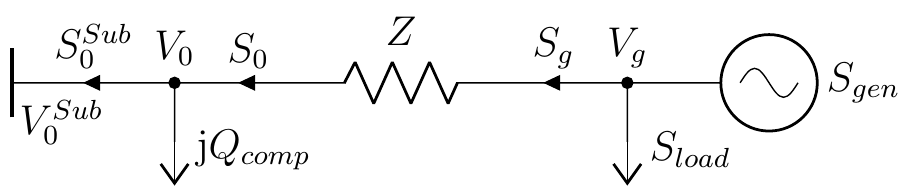}
\caption{Model accounting for a distribution system load and substation reactive power compensation.}\label{f:pscc_gen_model}
\end{figure}

In order to characterise the impact on a real distribution system, we consider the network shown in Figure \ref{f:pscc_gen_model}. The preceding analysis holds, with 
\begin{align*}
S_{0}^{Sub} =& \, S_{0} - \text{j}Q_{comp}\,,\\
S_{g} =& \, S_{gen} - S_{load}\,.
\end{align*}
Here, $S_{load}$ is the total feeder load at rated voltage, and $Q_{comp}$ is the compensation reactive power to improve the power factor of the power transferred through the substation $S_{0}^{Sub}$. As such, we append an additional constraint to \eqref{e:main_1} such that we also do not overload the substation transformer, i.e. 
\begin{equation}
P_{0}^{Sub} \leq P_{+}\,.
\end{equation}

\subsection{Three phase calculations of $Z$}

To calculate an approximate two bus equivalent for $Z$ for a general network, we utilise the properties of the impedance matrix $Z_{bus} \in \mathbb{C}^{n\times n}$, which is defined by 
\begin{equation}\label{e:zbus}
V = Z_{bus}I\,,
\end{equation}
where $V,I \in \mathbb{C}^{n}$ are the node voltages and nodal current injections respectively. The Thevenin impedance between two nodes in single phase equivalent circuits is calculated by injecting currents $I_{bus} = \delta I$ at the buses in question, and using the calculated voltage drop $\delta V = V_{bus,2} - V_{bus,1}$ from \eqref{e:zbus}, then calculate $Z = \delta V/\delta I$ \cite{grainger1994power}. In the case of an unbalanced distribution network, we use a similar method, with positive sequence current $\delta I$ injected, and positive sequence voltage $\delta V$ measured, such that $Z$ remains a scalar quantity.

\section{Case Study}

\begin{figure}
\includegraphics[width=0.48\textwidth]{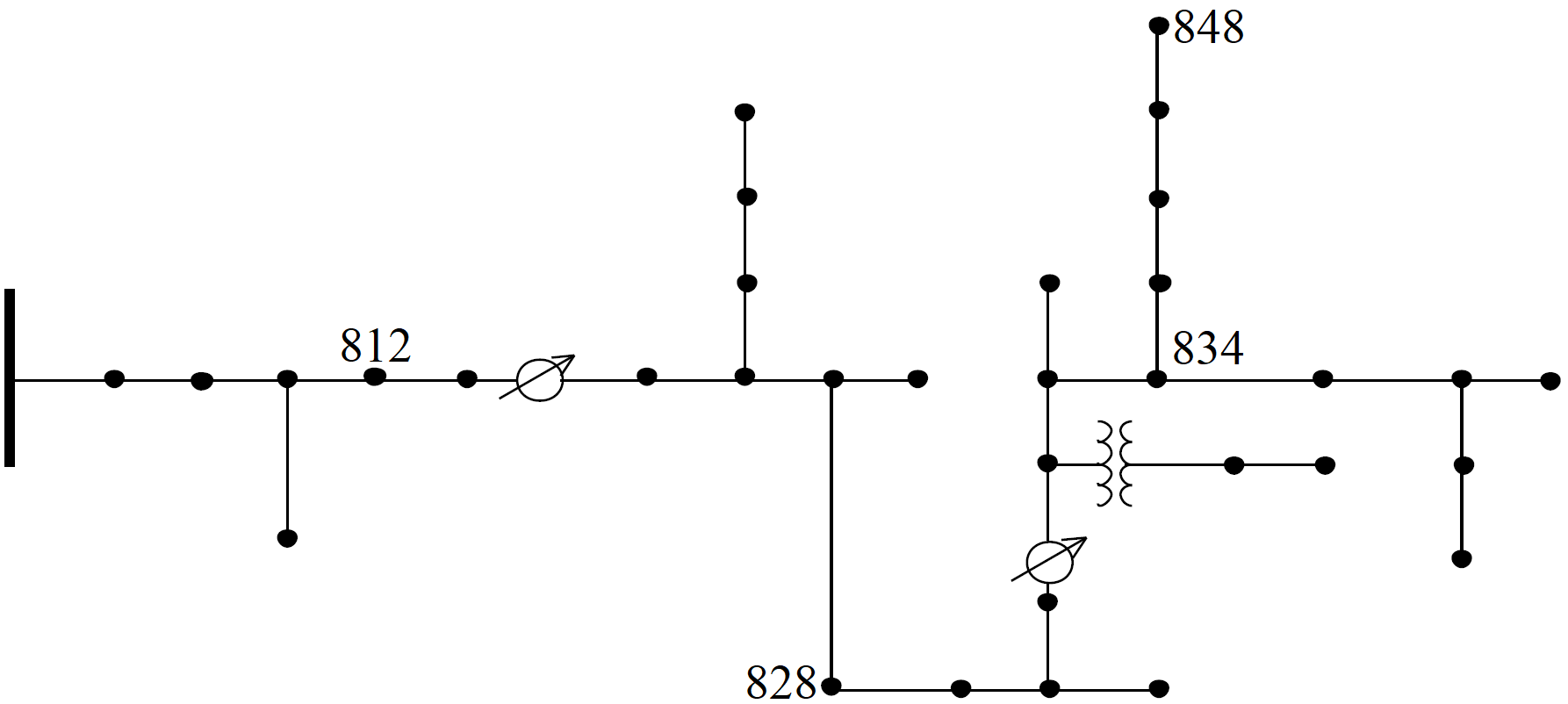}
\caption{IEEE 34 bus test distribution feeder and buses considered. Note that the voltage regulators shown are fixed for these results.}\label{f:ieee34_highlight}
\end{figure}

In order to evaluate the accuracy of the preceding analysis, we consider the network behaviour when a medium sized PV generator is placed at one of four buses in the IEEE 34 bus distribution test feeder (see Figure \ref{f:ieee34_highlight}) \cite{ieee2017distribution}. The feeder is modelled and solved in OpenDSS \cite{opendss2017}. The parameters calculated for the buses studied are given in Table \ref{t:cases}. The MATLAB/Octave and OpenDSS code used in this paper is available at:
\begin{center}
\texttt{https://github.com/deakinmt/pscc18}
\end{center}

\begin{table}[!t]
\caption{Network Parameters (pu). $S_{base} = 2.5$ MVA, $V_{base,LL} = 69$ kV.}
\label{t:cases}
\centering
\begin{tabular}{l l l l l l l l}
\toprule
Bus & $\lambda$  & $|Z|$ & $V_{0}$ & $V_{+}$ & $I_{+}$ (A) & $P_{+}$ & $S_{load}$ \\
\midrule
812 & 1.41 & 0.078 & 1.05 & 1.06 & 180 & 1.0 & $0.72\angle 9.4^{\circ}$\\
828 & 1.52  & 0.119 & 1.05 & 1.06 	& 180 & 1.0 & $0.72\angle 9.4^{\circ}$\\
834 & 1.85 & 0.203 & 1.05 & 1.06 	& 180 & 1.0 & $0.72\angle 9.4^{\circ}$\\
848 & 1.87 & 0.212 & 1.05 & 1.06 & 180 & 1.0 & $0.72\angle 9.4^{\circ}$ \\
\bottomrule
\end{tabular}
\end{table}

\subsection{Detailed study: Bus 834}\label{s:detailed_study}

\begin{figure}
\centering
\subfloat[]{\includegraphics[width=0.31\textwidth]{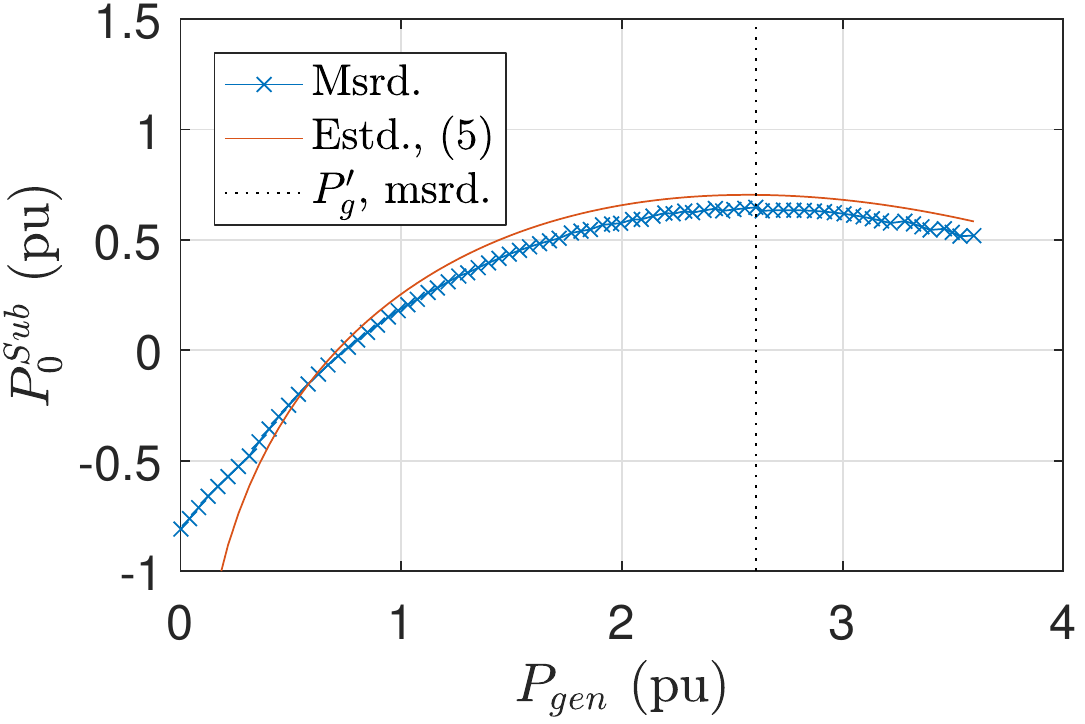}\label{f:34bus_pgp0}}\hfill
\subfloat[]{\includegraphics[width=0.31\textwidth]{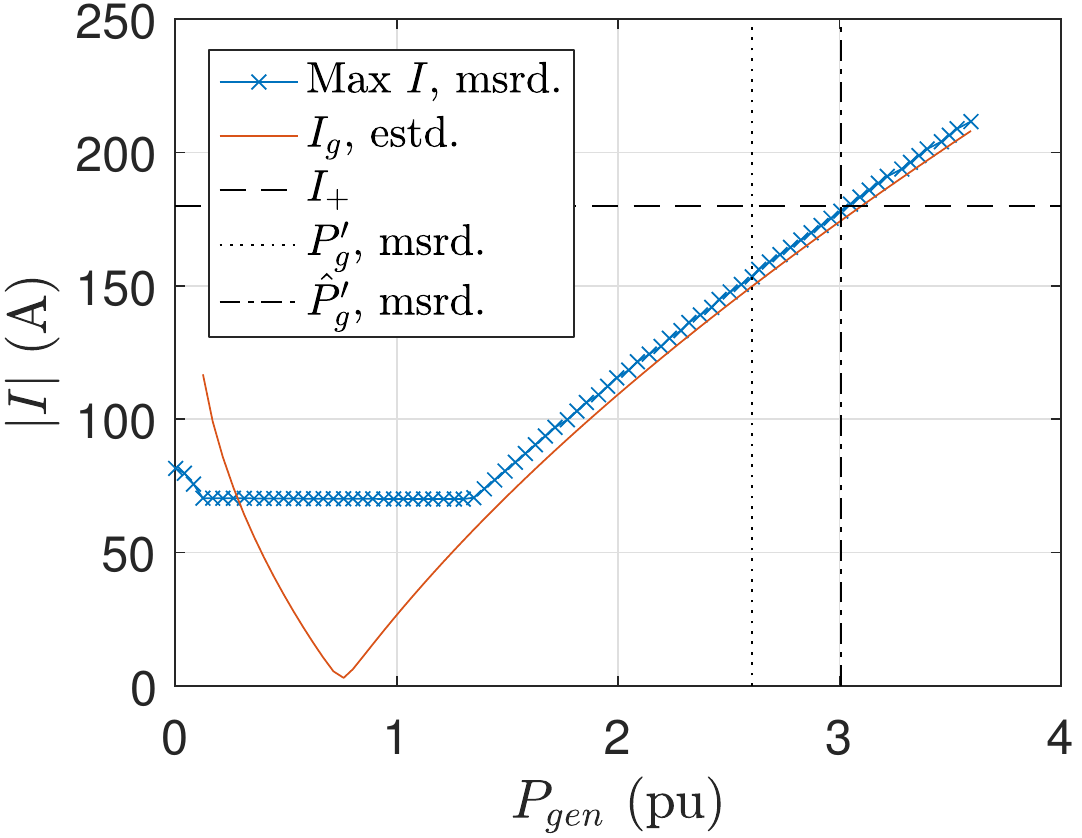}\label{f:34bus_imax}}\hfill
\subfloat[]{\includegraphics[width=0.31\textwidth]{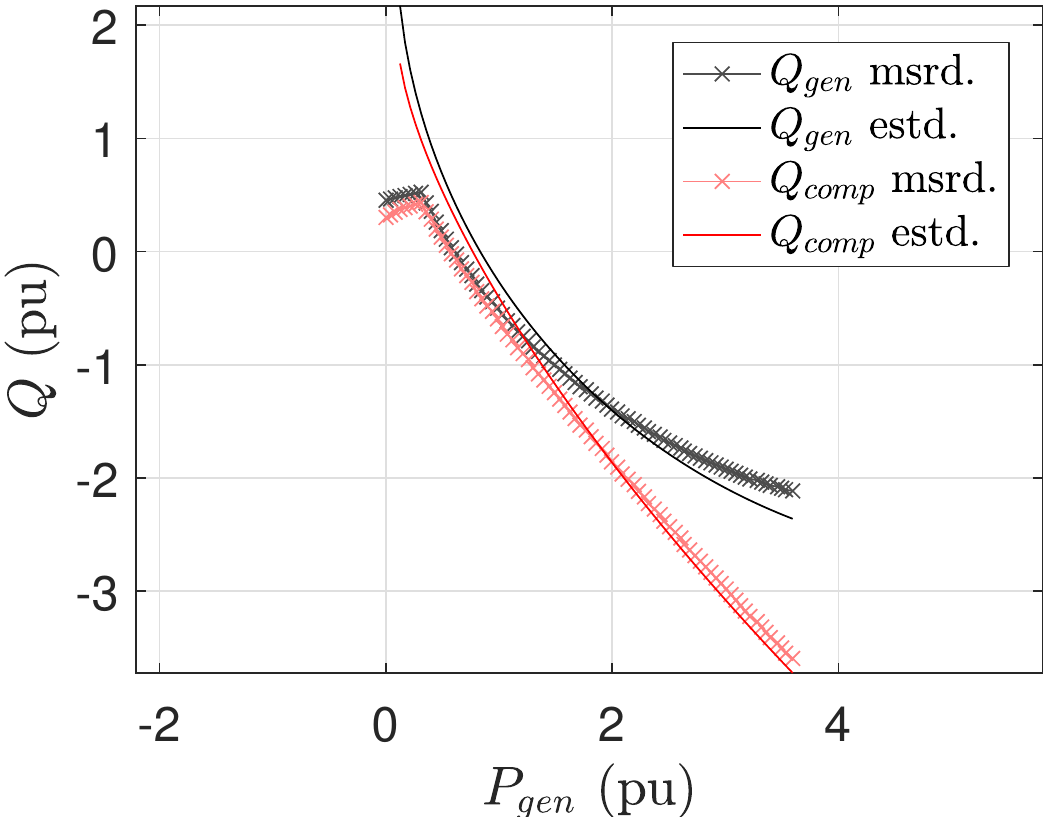}\label{f:34bus_pgqgq0}}
\caption{Feasible real and reactive power flows at bus 834 versus generated power $P_{gen}$. (a) $P_{gen}$ versus network transferred power $P_{0}^{Sub}$. (b) $P_{gen}$ versus maximum circuit current (estimated current defined by \eqref{e:eqn0}). (c) $P_{gen}$ versus estimated and measured reactive powers (estimated values defined by \eqref{e:v2}.)}
\end{figure}

We first demonstrate the existence of the loss induced power transfer limits for a bus close to the end of the feeder. To do so, a wide range of real and reactive powers are generated at the bus in question, and any infeasible points (i.e. any points at which upper voltage limits are violated) are removed. For each value of real power generated, the reactive power is chosen which maximises the real power transferred through the substation. This procedure yields the curve shown in Figure \ref{f:34bus_pgp0}. In addition, the solution to the fixed P-Q curve (\eqref{e:v2} at fixed $|V_{g}|$) is also shown on this figure for the network parameters given. We see that this corresponds well to the behaviour seen at high generator powers, and that the constraint on $P_{+}$ is not violated for any $P_{gen}$. Low generator powers see a divergence from this behaviour, as the network is not operating at the voltage constraint at these points.

Figure \ref{f:34bus_imax} shows the maximum network phase current and predicted currents (with permissible phase currents for conductors taken from \cite{kersting2001radial}). Vertical lines indicate the measured values of $P_{g}'$ and $\hat{P}_{g}$. We see that the network is not overloaded at the measured value of $P_{g}'$, as the maximum network current is below the permitted value. Therefore, we can conclude that the problem posed by \eqref{e:main_1} is indeed given by the marginal loss-induced limit $P_{g}'$, and that operating at the thermal limit $\hat{P}_{g}$ would yield less power transfer due to increased losses. Finally, in Figure \ref{f:34bus_pgqgq0}, we see that the estimated and measured reactive power flows are estimated reasonably well across a range of generated powers (again diverging at low generation).

\subsection{Model Accuracy}

\begin{figure}
\centering
\includegraphics[width=0.44\textwidth]{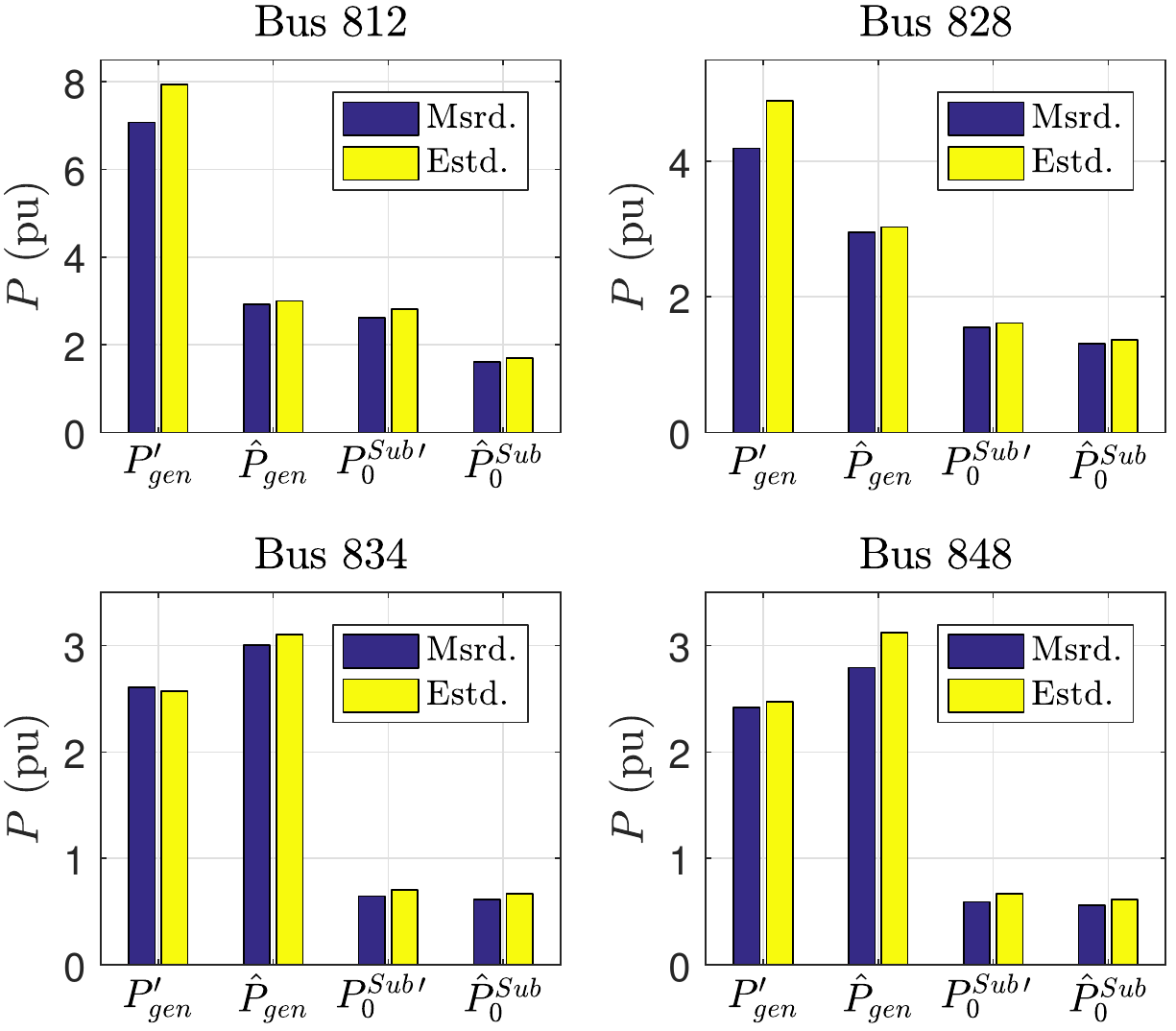}
\caption{The marginal and thermal induced power transfer limits for four buses of the IEEE 34 bus network. $S_{base} = 2.5$ MVA.}\label{f:34bus_accuracy}
\end{figure}

To consider the accuracy of the model, we consider the calculation of $P_{g}'$ and $\hat{P}_{g}$ for a generator located at a range of buses in the network, repeating the analysis of Section \ref{s:detailed_study}. The results are plotted in Figure \ref{f:34bus_accuracy}, and the errors between predicted and measured values shown in Table \ref{t:res_error}. We see that the theorem is able to accurately predict the solution of the problem to a relatively good accuracy in all buses. 

In particular, we note that the qualitative nature of the results is accurate using the two-bus model. There will always be a need for detailed simulations; however, we have shown that to a good degree of accuracy the power flows can be estimated well with the theory presented here.

\begin{table}                                                                                                                       
\caption{Marginal and thermal power limit error (pu)}                                                                               
\label{t:res_error}                                                                                                             \centering                                                                                                                          
\begin{tabular}{lllll}                                                                                                              
\toprule                                                                                                                            
Bus & $\epsilon (P_{gen}')$ & $\epsilon (\hat{P}_{gen})$ & $\epsilon (P_{0}^{Sub \, \prime})$  & $\epsilon (\hat{P}_{0}^{Sub})$   \\
\midrule                                                                                                                            
812 & -0.86 & -0.07 & -0.19 & -0.09 \\                                                                                              
828 & -0.70 & -0.08 & -0.06 & -0.06 \\                                                                                              
834 & 0.03 & -0.10 & -0.06 & -0.05 \\                                                                                               
848 & -0.06 & -0.33 & -0.07 & -0.05 \\                                                                                              
\bottomrule                                                                                                                         
\end{tabular}                                                                                                                       
\end{table}

\section{Conclusions}

The maximum power transfer theorem, the voltage `stability' boundary, and unity power factor control have all been demonstrated to be inadequate in the calculation of a general maximum power transfer limit for generation in distribution networks. A theorem has therefore been presented to solve this problem, subject to voltage and thermal constraints. The solution is in closed form and parametrised in terms of the $R/X$ ratio.

In addition to the familiar `thermal' loss induced limits, we demonstrate for the first time the existence of a `marginal' loss induced power transfer limit; that is, the point at which losses increase at a faster rate than the generated power. This bounds the amount of zero marginal cost energy that should be generated on a feeder, and as such, a point at which generated power should always be curtailed. On the other hand, it also clearly demonstrates how the `cost' of curtailment varies with the amount of power generated (and, indeed, is negative as we cross the marginal power transfer limit). The method also has the advantage of giving bounds on the maximum reactive power that would ever be required to manage these reactive power flows. Finally, we have demonstrated that, if it is possible to design a suitable control scheme, that it may be advantageous to operate in the `low voltage/high current' region for networks with small $R/X$ ratios.

These results are studied on the unbalanced three-phase IEEE 34 bus distribution feeder and the two cases are shown to exist, with accurate results on buses across the network. In practise, there are a wide range of technical constraints that must be managed when considering distribution networks. This paper shows that accurate, closed-form solutions to certain power flow problems do exist and can predict optimal network behaviour accurately. These types of solutions have the advantage of being extremely fast and allowing mathematical methods such as calculus to be used to analyse network behaviour.

\section*{Appendix: Proof of Theorem \ref{th:main_result}}\label{s:app1}

First we note that by expanding \eqref{e:v2} we obtain
\begin{equation}
|V_{g}|^{4} - (V_{0}^{2} + 2\tilde{P}_{g})|V_{g}|^{2} + \tilde{P}_{g}^{2} + \tilde{Q}_{g}^{2} = 0 \,, \label{e:id_2}
\end{equation}
\noindent irrespective of the operating region. From \eqref{e:rot_scale} and \eqref{e:id_1},
\begin{align}
P_{0} 	=& \, \dfrac{1}{|Z|^{2}}(R\tilde{P}_{0} - X\tilde{Q}_{0}),\\
		=& \, \dfrac{1}{|Z|^{2}}\bigg((R(|V_{g}|^{2} - V_{0}^{2}) - R\tilde{P}_{g} - X\tilde{Q}_{g}\bigg)\,.\label{e:pgtqgt_signdrn}
\end{align}
\noindent We assume for now that the maximum power that can be transferred is found on the line $|V_{g}| = V_{+}$. At function extrema we can differentiate such that
\begin{align*}
\dfrac{dP_{0}}{d\tilde{P}_{g}} =& \, \lambda + \dfrac{d\tilde{Q}_{g}}{d\tilde{P}_{g}}\,,\\
			=& \, \lambda + \dfrac{V_{+}^{2} - \tilde{P}_{g}}{\tilde{Q}_{g}}\,,\\
			=& \, 0\,, \\
			\Rightarrow \lambda \tilde{Q}_{g} =& \, \tilde{P}_{g} - V_{+}^{2} \, ,\\
			\Rightarrow \lambda ^{2}\tilde{Q}_{g}^{2} =& \, \tilde{P}_{g}^{2} - 2 \tilde{P}_{g}|V_{g}|^{2} + |V_{g}|^{4}.
\end{align*}
\noindent Using \eqref{e:id_2} we can therefore write
\begin{equation*}
\lambda ^{2} \bigg( (V_{0}^{2} + 2 \tilde{P}_{g})V_{+}^{2} - \tilde{P}_{g}^{2} - V_{+}^{4} \bigg) = - 2 \tilde{P}_{g}V_{+}^{2} + \tilde{P}_{g}^{2} + V_{+}^{4} \, .
\end{equation*}
\noindent This is quadratic in $\tilde{P}_{g}$ with solutions
\begin{equation*}
\tilde{P}_{g} = V_{0}V_{+}(\dfrac{V_{+}}{V_{0}} \pm \dfrac{R}{|Z|})\,. \label{e:pg_opt_quadratic}
\end{equation*}
We use this result with \eqref{e:id_2} to show
\begin{equation*}
\tilde{Q}_{g} = \pm V_{0}V_{+}\dfrac{X}{|Z|}\,.
\end{equation*}

\noindent By \eqref{e:pgtqgt_signdrn}, we take (-) instead of (+) for both $\tilde{P}_{g}$ and $\tilde{Q}_{g}$ to maximise the net power sent to the grid, yielding $(\tilde{P}_{g}',\tilde{Q}_{g}')$.

Finally, we must show that \eqref{e:PVlocus} holds. From \eqref{e:pgtqgt_signdrn}, we can substitute back in $\tilde{P}_{g}',\tilde{Q}_{g}'$ to show
\begin{equation*}
P_{0}' = \dfrac{V_{0}^{2}}{Z}(\dfrac{V_{+}}{V_{0}}- \dfrac{R}{|Z|})\,.
\end{equation*}
\noindent Therefore, any increase in $|V_{g}| = V_{+}$ increases the maximum power transfer.

\section*{Acknowledgements} We wish to thank Dr. Frank Johnson for motivating the exploration of themes in this paper, and to the Sir John Aird and Clarendon Scholarships for their generous support.

\bibliographystyle{IEEEtran}
\bibliography{pscc18_bib}{}


%\begin{thebibliography}{1}
%\bibitem{Shell}
%M.~Shell, \emph{How to Use the IEEEtran Latex Class}, Latex Archive Contents, \verb+http://www.ieee.org/conferences_events/+ \verb+conferences/publishing/templates.htm+
%\end{thebibliography}

%\bibitem{Shell2}
%M.~Shell, \emph{How to Use the IEEEtran Latex Class}, Latex Archive Contents, \verb+http://www.ieee.org/conferences_events/+ \verb+conferences/publishing/templates.htm+
%
%\bibitem{IEEEhowto:kopka2}
%\end{thebibliography}

\end{document}